\documentstyle[12pt]{article}

\def\R{{I\!\!R}}

\def\CC{{\rm \kern.24em \vrule width.02em height1.4ex
    depth-.05ex \kern-.26em C}}

\def\TagOnRight

\def\AA{{\it I}\hskip-3pt{\tt A}}

\def\QQ{\rlap {\raise 0.4ex \hbox{$\scriptscriptstyle |$}}
  {\hskip -0.1em Q}}
\def\RR{\rm I\hskip-1.8pt R}
\def\ZZ{\mbox{{\rm Z\hskip-4pt \rm Z}}}   %\def\ZZ{Z\!\!\!Z}

\newcommand{\be}{\begin{equation}}
\newcommand{\ee}{\end{equation}}
\newcommand{\bea}{\begin{eqnarray}}
\newcommand{\eea}{\end{eqnarray}}
\newcommand{\Bea}{\begin{eqnarray*}}
\newcommand{\Eea}{\end{eqnarray*}}

\catcode`\@=11
\def\theequation{\@arabic{\c@equation}}
\catcode`\@=12

\newcommand{\bi}{\begin{itemize}}
\newcommand{\ei}{\end{itemize}}

\newtheorem{Definition}{Definition}[section]
\newtheorem{Theorem}[Definition]{Theorem}

\newtheorem{Proposition}[Definition]{Proposition}

\renewcommand{\theequation}{ \thesection.\arabic{equation} }

\newcommand{\C}{\mbox{\rule{.1mm}{2.5mm}$\!\!C$}}%

\begin{document}
\parskip 12pt
\parindent 0pt
\renewcommand{\thepage}{\arabic{page}}

\begin{center}

{\large {\bf Probabilistic Representations of Solutions \\
of the Forward Equations} \\
\vskip 1em By}
\end{center}
\noindent
\begin{tabular}{ll}
  B. Rajeev,  & ~~~~~~~~~S. Thangavelu, \\
  Statistics and Mathematics Unit, & ~~~~~~~~~Department of Mathematics,\\
  Indian Statistical Institute, & ~~~~~~~~~Indian Institute of Science, \\
    Bangalore - 560 059. & ~~~~~~~~~Bangalore - 560 012. \\
  brajeev@isibang.ac.in & ~~~~~~~~~veluma@math.iisc.ernet.in \\
  \end{tabular}

{\bf Abstract:} In this paper we prove a stochastic representation
for solutions of the evolution equation \Bea
\partial_t \psi_t &=& \frac{1}{2} L^\ast \psi_t \Eea where
$L^\ast$ is the formal adjoint of a second order elliptic
differential operator $L$, with smooth coefficients, corresponding
to the infinitesimal generator of a finite dimensional diffusion
$(X_t)$. Given $\psi_0 = \psi $, a distribution with compact
support, this representation has the form $\psi_t = E(Y_t(\psi))$
where the process $(Y_t(\psi))$ is the solution of a stochastic
partial differential equation connected with the stochastic
differential equation for $(X_t)$ via Ito's formula.


\section{{\protect\large {\bf Introduction }}}
\setcounter{equation}{0}

The first motivation for the results of this paper is that they
extend the results in \cite{RT} for Brownian motion, to more general
diffusions. To recall, there we had recast the classical
relationship between Brownian motion and the heat equation in the
language of distribution theory. To be more precise, the solutions
of the initial value problem, \Bea
\partial_t \psi_t &=& \frac{1}{2} \Delta \psi_t,  \\
\psi_0 &=& \psi \Eea for any distribution $\psi \in {\cal S}'
(\R^d)$ were represented in terms of a standard $d$-dimensional
Brownian motion $(X_t)_{t \geq 0}$, as $\psi_t =E\tau_{X_t}
(\psi)$. Here $\tau_x :\R^d \rightarrow \R^d$ is the translation
operator and the expectation is taken in a Hilbert space $S_p
\subset {\cal S}'$ in which the process takes values.

The main result of this paper is that the above result extends to
the solutions of the initial value problem \bea
\partial_t \psi_t &=& \frac{1}{2} L^\ast \psi_t, \nonumber \\
\psi_0 &=& \psi \eea where $\psi \in {\cal E}'$, i.e it is a
distribution on $\R^d$ with compact support and $L^\ast$ is the
formal adjoint of $L$, a second order elliptic differential operator
with smooth coefficients given as the first component of the pair
$(L,A)$, where \Bea L &=& \frac{1}{2} \sum\limits_{i,j}
(\sigma\sigma^t )_{ij} (x) ~\partial^2_{ij}
+\sum\limits_i b_i (x)~\partial_i,\\ A &=& (A_1,\cdots , A_d), \\
A_k &=& - \sum\limits_{i=1}^d \sigma_{ij}(x)
\partial_i. \\ \Eea  To describe the
stochastic representation of solutions of (1.1), let $(X(t,x))$
denote the solutions of the stochastic differential equation \Bea
dX_t &=& \sigma (X_t)\cdot dB_t
+ b(X_t) ~dt ,\\
X_0 &=& x. \Eea  Then it is well known that a.s, $x \rightarrow
X(t,x) $ is smooth and induces a map $X_t : C^\infty \rightarrow
C^\infty $, namely $X_t(\phi)(x) = \phi(X(t,x))$ for $\phi \in
C^\infty$. Here and in the rest of the paper $C^\infty$ denotes the
space of all smooth functions on $\R^d$. Let $Y_t : =X_t^\ast : {\cal E}'
\rightarrow {\cal E}'$ be the adjoint of the map $X_t : C^\infty
\rightarrow C^\infty$. If $\psi \in {\cal E}' \subset {\cal S}'$
where ${\cal S}'$ is the space of tempered distributions, then we
can show that the process $(Y_t(\psi))$ takes values in one of the
Hilbert spaces $S_{-p}$,$ p > 0$ that define the countable 
Hilbertian structure of ${\cal S}'$
 (Proposition 3.1). Our stochastic representation now reads,
$\psi_t=E(Y_t (\psi))$(Theorem 4.3). These results  extend the well
known results connecting diffusion processes and PDE (see
\cite{B}, \cite{D}, \cite{SV}, \cite{F}, \cite{BL}). Moreover, they
establish a natural link with the subject of stochastic partial
differential equations (see \cite{K},\cite{I}, \cite{KR},\cite{KMW}
\cite{W},\cite{DZ}), viz. the process $(Y_t (\psi))$ for $\psi \in
{\cal E}'$, is the solution of a stochastic partial differential
equation (eqn (3.7) below) associated naturally with the equation
for $(X_t)$ via Ito's formula. This  stochastic partial differential
equation  is different from the one satisfied by the process
$(\delta_{X_t})$ in \cite{R},\cite{GMR},\cite{GMRMon} - the former
is associated with the operators $(L^\ast, A^\ast)$ as above, the
latter being associated with the random  operators
$(L(t,\omega),A(t,\omega))$, $A(t,\omega) := (A_1(t,\omega),\cdots
A_d(t,\omega))$,   \Bea A_k(t, \omega) &=&  - \sum\limits_{i=1}^d
\sigma_{ij}(X_t(\omega)
\partial_i, \\  L(t,\omega) &=& \frac{1}{2} \sum\limits_{i,j}
(\sigma\sigma^t )_{ij} (X_t(\omega)) ~\partial^2_{ij} -\sum\limits_i
b_i (X_t(\omega))~\partial_i    \Eea (see \cite{R},\cite{GMRMon}).
However it is easily seen that when $\psi=\delta_x$, the process
$(Y_t (\psi))$ is the same as the process $(\delta_{X_t})$. We also
note that solutions $(\psi_t)$ of equation (1.1) are obtained by
averaging out the diffusion term in the stochastic partial
differential equation  satisfied by $(Y_t(\psi))$ (see Theorem 4.3),
a result that corresponds quite well with the original motivation
for studying stochastic partial differential equations, viz.
`Stochastic PDE = PDE + noise' (see \cite{W} for example).

The definition of $ Y_t $ as the adjoint of the map $X_t : C^\infty
\rightarrow C^\infty$ induced by the flow $(X(t,x,\omega))$ does not 
automatically lead to good path properties
for the process $(Y_t(\psi)), \psi \in {\cal E}'$. To get these we
generalise the representation \Bea Y_t (\psi) (\omega) = \int
\delta_{X(t,x,\omega)} ~d \psi (x) \Eea which is easily verified
when $\psi$ is a measure with compact support. Indeed, \Bea \big
\langle Y_t (\psi) (\omega), \varphi \rangle &=&
\langle \psi, X_t (\omega) \circ \varphi \rangle \\
&=& \int X_t (\omega) \circ \varphi (x)~d \psi (x)  \\
&=& \int \varphi (X(t,x,\omega)) ~d \psi (x) \\
&=& \int \langle  \delta_{X(t,x,\omega)},\varphi \rangle ~d \psi (x)
\\
&=& \langle  \int \delta_{X(t,x,\omega)} ~d\psi (x),\varphi \rangle.
\Eea Here the integral $\int \delta_{X(t,x,\omega)} d\psi(x)$ is
understood in the sense of Bochner and takes values in a suitable
Hilbert space $S_{-p} \subseteq {\cal S}'$. In Section 3, we define
the process $(Y_t(\psi))$, for $\psi \in {\cal E}'$, via a
representation such as the one above and verify that indeed $Y_t
(\psi)=X_t^\ast (\psi)$. In Theorem (3.3) we show that $(Y_t(\psi))$
satisfies a stochastic partial differential equation.

In Section 2, we state some well known results from the theory of
stochastic flows \cite{K} in a form convenient for our purposes.
These are used in Section 3 for constructing the process $(Y_t
(\psi))$ and proving its properties. In Section 4, we prove the
representation result for solutions of equation (1.1). Our results
amount to a proof of existence of solutions for equation (1.1). Our
proofs require that the coefficients be smooth. However, we do not
require that the diffusion matrix be non degenerate. The uniqueness
of solutions of (1.1) can be shown to hold if the so called
`monotonicity inequality' (see (4.2) below) for the pair of
operators $(L^\ast,A^\ast)$ is satisfied (Theorem 4.4). It may be
mentioned here that the `monotonicity inequality' is known to hold
even when the diffusion matrix is degenerate (see \cite{GMRMon}). As
mentioned above, when $\psi=\delta_x, x \in \R^d, Y_t
(\psi)=\delta_{X_t}$ where $X_0=x$. As discussed above, for $\psi
\in {\cal E}'$, the solutions $(Y_t(\psi))$ of the stochastic partial 
differential equation (3.7) can be constructed out of the
particular solutions $(\delta_{X_t})$ corresponding to $\psi
=\delta_x$. In other words the processes $(\delta_{X_t}), X_0 = x$,
can be regarded as the `fundamental solution' of the stochastic
partial differential equation (3.7) that $(Y_t (\psi))$ satisfies.
This property is preserved on taking expectations: in other words
$E\delta_{X_t}, X_0 = x$ is the fundamental solution of the  partial
differential equation  (1.1) satisfied by $\psi_t = EY_t(\psi)$ - a
well known result for probabilists, if one notes that
$E\delta_{X_t} = P(t,x,\cdot)$ the transition probability measure of
the diffusion $(X_t)$ starting at $x$ (Theorem 4.5). Assuming that
this diffusion has a density $p(t,x,y)$ satisfying some mild
integrability conditions, we deduce some well known results. If
$L^\ast =L$ then the density is symmetric (Theorem (4.6)), and in
the constant coefficient case we further have $p(t,x,y)
=p(t,0,y-x)$. Finally if $T_t : C^\infty \rightarrow C^\infty$
denotes the semigroup corresponding to the diffusion $(X_t)$, and
$S_t : {\cal E}'\rightarrow {\cal E}'$ is the adjoint, then $S_t$,
given by $S_t (\psi)=EY_t(\psi)$, is a uniformly bounded (in $t$)
operator when restricted to the Hilbert spaces $S_{-p}$ (Theorem
(4.8)).

\section{{\protect\large {\bf Stochastic Flows}}}
\setcounter{equation}{0}

Let  $\Omega = C([0,\infty),\R^r)$ be the set of continuous functions
on $[0,\infty)$ with values in $\R^r$. Let $ ~{\cal F} $ denote the
Borel $\sigma$-field on $\Omega$ and let $P$ denote the  Wiener
measure. We denote $B_t(\omega ):=\omega (t), ~\omega \in \Omega, ~t
\geq 0$  and recall that under $P$, $(B_t)$ is a standard $r $
dimensional Brownian motion. Let $(X_t)_{t \geq 0}$ be a strong
solution on $(\Omega ,{\cal F}, P) $ of the stochastic differential
equation  \bea \left.
\begin{array}{llll}
dX_t &=& \sigma (X_t) \cdot dB_t+b(X_t) dt  \\
X_0 &=& x \end{array} \right\} \eea with $\sigma =(\sigma_{j}^i)$,
$i=1 \ldots d, j = 1 \ldots r$ and $b=(b^1, \ldots b^d)$, where
$\sigma^i_j$ and $b^i$ are  given by $C^\infty$ functions on $\R^d$
with bounded derivatives satisfying \Bea \|\sigma (x)\|+\|b(x)\|=
\left( \sum\limits^d_{i=1} \sum\limits^r_{j=1} |\sigma _{j}^i (x)|^2
\right)^{1/2} + \left( \sum\limits^d_{i=1} |b^i(x)|^2 \right)^{1/2}
\leq K(1+|x|) \Eea for some $K >0$. Under the above assumptions on
$\sigma$ and $b$, it is well known that a unique, non-explosive
strong solution $(X(t,x,\omega ))_{t \geq 0, x \in \R^d}$ exists on
$(\Omega, {\cal F}, P)$ (see \cite{IW}). We also have the following
theorem (see \cite{JB}, \cite{K} and \cite{IW}, p.251  ).

\begin{Theorem} For $x \in \R^d$ and $t \geq 0$, let $ (X(t,x,\omega ))$
be the unique strong solution of equation (2.1)  above. Then there exists
a process $(\tilde X (t,x,\omega ))_{t \geq 0, x \in \R^d}$  such that\\
(1) For all $ x \in \R^d, P\{\tilde X (t,x,\omega )=X(t,x,\omega ),~\forall~
t \geq 0\}=1$. \\
(2) For a.e. $\omega (P), x \rightarrow \tilde X (t,x,\omega )$ is a
diffeomorphism for all $t \geq 0$.\\
(3) Let $\theta _t: \Omega \rightarrow \Omega$ be the shift operator
i.e. $\theta _t \omega (s)=\omega (s+t)$; then for $s,t \geq 0$, we have
$$\tilde X(t+s,x,\omega )=\tilde X(s,\tilde X (t,x,\omega ),\theta _t \omega ) $$
for all $x \in \R^d$, a.e. $\omega$ $(P)$.
\end{Theorem}

Denote by $\partial X(t, x ,\omega )$ the $d \times d$ matrix valued
process, given for a.e. $\omega $ by \linebreak $(\partial X)^i_j
(t,\omega ,x) =\frac{\partial X^i}{\partial x^j}(t, x, \omega )$ for
all $t \geq 0$ and $x \in \R^d$. Denote by $\partial \sigma _\alpha
(x)$ the $d \times d$ matrix $(\partial \sigma _{\alpha }
(x))^i_j=\frac{\partial \sigma _\alpha^i (x)} {\partial x^j}$ and by
$\partial b(x)$ the $d \times d$ matrix $(\partial b(x))^i_j=
\frac{\partial b^i (x)}{\partial x^j}.$ Then it is well known (see
\cite{K}) that almost surely, $\partial X(t,\omega ,x)$ is
invertible for all $t$ and $x$ and the inverse satisfies the SDE
\Bea dJ_t &=& - \sum\limits_{\alpha = 1}^r  J_t \cdot
\partial \sigma _\alpha (X_t) ~dB^\alpha
_t \\
&& - J_t \cdot \left[\partial b(X_t)-\sum\limits_{\alpha = 1}^r
(\partial \sigma _
\alpha ) \cdot (\partial \sigma _\alpha) (X_t) \right] dt, \\
J_0 &=& I \Eea where $I$ is the identity matrix and $J_t \cdot
\partial \sigma_\alpha (X_t)$ etc. denote the product of $d \times
d$ matrices. In proving our results, we will need to show that
$\sup\limits_{x \in K}|\partial^r X(t,x)|^q$ and $\sup\limits_{x \in
K} |(\partial X)^{-1} (t,x)|$ (here $ \partial^r :=
\partial^{r_1}_1 \ldots \partial^{r_d}_d$ and $|\cdot |$ denotes the
Euclidean norm on $\R^d$ in the first case and on $\R^{d^2}$ in the
second case) have finite expectation for $q \geq 1$ and $K \subseteq
\R^d$ a compact set. To do this we will use the results of section
4.6 of \cite{K}, as also the notation there. First, we note that the
stochastic differential equations for $(X_t)$ and $(\partial X(t))$
can be combined into a single stochastic differential equation  in
$\R^{d+d^2}$, which in the language of \cite{K}, can be based on a
spatial semi-martingale $F(x,t)= (F^1(x,t), \ldots F^{d+d^2}
(x,t))$. Having done this and having verified the regularity
hypothesis on the local characteristics of $F(x,t)$ we can apply
Corollary 4.6.7 of \cite{K} to get our results.

To fix notation we note that the set $\{k :d < k \leq d^2+d\}$ is in
1-1 correspondence with $\{(i,j):1 \leq i, j \leq d\}$. We fix such
a correspondence and write $k \leftrightarrow (i,j)$ for $d < k \leq
d+d^2$ and $1 \leq i, j \leq d$. If $x \in \R^{d+d^2}$ we will write
$x=(x_1,x_2)$ where $x_1 \in \R^d$ and $x_2 \in \R^{d^2}$. For $1 \leq k \leq 
d$, let $$F^k(x,t) =\sum\limits^r_{\alpha =1} \sigma ^k_\alpha (x_1) 
B_t^\alpha + b^k(x_1) t.$$ For $d +1 \leq k \leq d+d^2$, let 
\Bea F^k (x,t) &=&
-\sum\limits_{\alpha = 1}^r (x_2 \cdot \partial \sigma_\alpha
(x_1))^i_j ~B^\alpha _t\\
&-& t \left(x_2 \cdot \left[\partial b(x_1)-\sum\limits_{\alpha =
1}^r (\partial \sigma _\alpha ) \cdot (\partial \sigma _\alpha)(x_1)
\right] \right)^i_j \Eea where $k \leftrightarrow (i,j), x_2 \cdot
\partial \sigma _\alpha (x_1)$ is the product of $d \times d$
matrices $x_2$ and $\partial \sigma _x (x_1)$ etc. The local
characteristics of $F(x,t)$ are then given by $(\alpha (x,y,t),
\beta (x,t),t)$ where: \Bea
\beta ^k (x,t) &=& b^k (x_1), ~~1 \leq k \leq d \\
&=& x_2 \cdot \partial b(x_1) -\sum\limits_{\alpha = 1}^r (\partial
\sigma _\alpha)
\cdot (\partial \sigma _\alpha ) (x_1), ~~ d+1 \leq k \leq d+d^2\\
\beta (x,t) &=& (\beta ^1 (x,t), \ldots \beta ^{d+d^2} (x,t)),
~~x=(x_1,x_2) \in \R^{d+d^2}. \Eea Further, for $ x = (x_1 , x_2 ), y
= (y_1 , y_2 )$ ,$\alpha (x,y,t)=\alpha ^{k \ell} (x,y,t), ~1\leq k,
~\ell \leq d+d^2$ where \Bea \alpha ^{k \ell} (x,y,t) &=& (\sigma
(x_1)\cdot \sigma^t (y_1))^k_\ell,
~~ 1 \leq k, \ell \leq d \\ \\
&=& \sum\limits_{\alpha = 1}^r (x_2 \cdot \partial \sigma _\alpha
(x_1))^i_j
(y_2 \cdot \partial \sigma _\alpha (y_1))^{i'}_{j'} \\
&& d+1 \leq k,
\ell \leq d+d^2, k \leftrightarrow (i,j), \ell\leftrightarrow (i',j') \\
&& \\
 &=& -\sum\limits_{\alpha =1}^r \sigma_\alpha ^k  (x_1) (y_2 \cdot
\partial
\sigma _\alpha  (y_1))^i_j \\
&& 1 \leq k \leq d, d+1 \leq \ell \leq d+d^2 \mbox{~and}~
\ell \leftrightarrow (i,j) \\
&&\\
 &=& -\sum\limits_{\alpha = 1}^r (x_2 \cdot \partial \sigma
_\alpha(x_1))^i_j
\sigma_\alpha ^\ell (y_1) \\
&& 1 \leq \ell \leq d, d+1 \leq k \leq d+d^2 \mbox{~and}~ k
\leftrightarrow (i,j). \Eea Let $n=d+d^2$ and $f: \R^n \rightarrow
\R^n, g :\R^n \times \R^n \rightarrow \R^{n^2}$. Consider for $m
\geq 1$ and $\delta >0$, the following semi-norms (see \cite{K}), 
\Bea \|f\|_{m,\delta } &=& \sup\limits_{x \in \R^n}
\frac{|f(x)|}{(1+|x|)} +\sum\limits_{1 \leq |\alpha |\leq m}
~~\sup\limits_{x \in \R^n} |\partial^
\alpha f(x)| \\
&& + \sum\limits_{|\alpha |=m} ~~\sup\limits_{\stackrel{x,y \in
\R^n}{x \neq y}}
\frac{|\partial^\alpha f(x) -\partial^\alpha f(y)|}{|x-y|^{\delta}} \\
\|g\|^\sim_{m,\delta} &=& \sup\limits_{x \in \R^n} \frac{|g(x,y)|}{(1+|x|)(1+|y|)}
+\sum\limits_{1 \leq |\alpha| \leq m} ~~\sup\limits_{x,y \in \R^n} |\partial^\alpha _
x \partial_y^\alpha g(x,y)|
\Eea
$$ + \sum\limits_{|\alpha |=m} ~~\sup\limits_{\stackrel{x,x',y,y'\in \R^n}{x
\neq x', y \neq y'}} \frac{|\partial_x^\alpha \partial^\alpha _y g(x,y)-
\partial^\alpha _x \partial^\alpha _y g(x',y)-\partial^\alpha _x \partial^\alpha _y
g(x,y')+\partial_x^\alpha \partial^\alpha _y g(x',y')|}{|x-x'|^\delta|y-y'|^\delta}.$$
Let for $0 < s<T$, $\{\varphi _{s,t} (x), s \leq t \leq T\}$ be the
solution of Ito's SDE based on the semi-martingale $F(x,t), x \in \R^{d+d^2}$, 
 i.e. $$\varphi _{s,t} (x) =x+\int\limits_s^{~~t} 
F(\varphi _{s,r} (x) ,dr). $$ Note that $\varphi ^k_{0,t}(x)=X^k(t,x_1)~~~1 
\leq  k\leq d$ and
$\varphi ^k_{0,t} (x)=(\partial X^{-1} )^i_j (t, x_1)$, $ d+1 \leq k
\leq d+d^2$ and $k\leftrightarrow (i,j), x =(x_1,x_2), x_1 \in \R^d$
and $x_2 \leftrightarrow I$ the identity matrix. We then have the
following theorem.

\begin{Theorem} Given $0 \leq s \leq T, \alpha =(\alpha _1 \ldots \alpha _d)$
a multi-index, $N >0$, and $q \geq 1 $, there exists
$C=C(s,T,\alpha ,N,q)>0$ such that
$$E \sup\limits_{|x|\leq N} |\partial^\alpha \varphi _{s,t} (x)|^q < C$$
for any $t$ satisfying $s \leq t \leq T$. In particular, for any
compact $K \subseteq \R^d, q \geq 1, \alpha $ a multi-index, there
exists $ ~ C>0$
$$ E \sup\limits_{x \in K} |\partial^\alpha X(t,x)|^q < C $$
and
$$ E \sup\limits_{x \in K} |\partial X^{-1} (t,x)^i_j|^q < C$$
for $0 \leq t \leq T$.
\end{Theorem}

{\bf Proof:} It is easily verified that the local characteristics $(\alpha ,
\beta ,t)$ verify for $m \geq 1, ~\delta =1$:
$$ \sup\limits_{t \leq T} \|\alpha (t)\|^\sim _{m,1} < \infty ~\mbox{and}~
\sup\limits_{t \leq T}\|\beta (t)\|_{m,1} < \infty.$$ In the
language of \cite{K}, the local characteristics $(\alpha ,\beta ,t)$
belong to the class $B^{m,1}_{ub}$ for all $m \geq 1$. Thus the
hypothesis of  Corollary 4.6.7 in \cite{K} is satisfied. Hence for
$p
>1, \alpha $ a multi-index, $N >0$ and $0 < s < T, ~\exists
~C=C(p,\alpha ,N,s,T)$
$$E \sup\limits_{|x| \leq N} |\partial^\alpha \varphi _{s,t} (x)|^{2p} < C.$$
 The result for $q=2p
>2$ and hence for $q \geq 1$ follows.
\hfill${\Box}$

\section{{\protect\large {\bf The Induced Flow on Distributions with compact
support}}}
\setcounter{equation}{0}

We will denote the modification obtained in Theorem 2.1 again by $(X
(t,x,\omega ))$. For $\omega$ outside a null set $\tilde{N}$, the
flow of diffeomorphisms induces, for each $t \geq 0$ a continuous
linear map, denoted by $X_t (\omega )$ on $C^\infty (\R^d)$. $X_t
(\omega ):C^\infty (\R^d) \rightarrow C^\infty (\R^d)$ is given by
$(X_t (\omega )(\varphi )(x)= \varphi (X(t,x,\omega ))$. This map is
linear and continuous with respect to the topology on $C^\infty
(\R^d)$ given by the following family of semi-norms: For $K
\subseteq \R^d$ a compact set, let $\|\varphi
\|_{n,K}=\max\limits_{|\alpha | \leq n}~~ \sup\limits_{x \in K}
|D^\alpha f(x)|$ where $\varphi \in C^\infty (\R^d)$ and $n \geq 1$
an integer and $\alpha =(\alpha _1, \ldots \alpha _d)$ and $|\alpha
|=\alpha _1 +\ldots +\alpha _d$. Let $K_{t,\omega }$ denote the
image of $K$ under the map $x \rightarrow X(t,x,\omega )$. Then
using the `chain rule' it is easy to see that there exists a
constant $C(t,\omega )
>0$ such that
$$\|X_t (\omega )(\varphi )\|_{n,K} \leq C(t,\omega )\|\varphi \|_{n,K_{t,\omega }}.$$
Let $X_t (\omega )^\ast$ denote the transpose of the map $X_t
(\omega ):C^ \infty (\R^d)\rightarrow C^\infty (\R^d)$. Then if
${\cal E}'$ denotes the space of distributions with compact support
we have $X_t (\omega )^\ast : {\cal E}'\rightarrow {\cal E}'$ is
given by
$$\langle X_t (\omega )^\ast  \psi, \varphi \rangle =\langle \psi, X_t (\omega )
\varphi \rangle$$ for all $\varphi \in C^\infty$ and $\psi \in {\cal
E}'$. Let $\psi \in {\cal E}'$. Let supp $\psi \subseteq K$ and let
$N=$ order$(\psi)+2d$. Then there exist continuous functions
$g_\alpha , |\alpha |\leq N, \mbox{supp}~ g_\alpha \subseteq V$
where $V$ is an open set having compact closure, containing $K$,
such that \bea \psi =\sum\limits_{|\alpha | \leq N} \partial^\alpha
g_\alpha . \eea See \cite{Tr}. Let $\varphi \in C^\infty (\R^d)$.
Let $f_i \in C^\infty (\R^d)$ and $f=(f_1, \ldots f_d)$. Then, using
the chain rule, it is easy to see that there exist polynomials
$P_\gamma (x_1 \ldots x_d)$, one for each multi-index $\gamma $,
with $|\gamma | \leq |\alpha |$ and $ {deg~} P_\gamma =|\gamma |$,
and such that \bea
\partial^\alpha (\varphi \circ f)(x)=\sum\limits_{\stackrel{|\gamma |\leq
|\alpha |}{|\beta |\leq|\alpha |} } (-1)^{|\gamma |} P_\gamma
(\partial^{\beta _1} f_1, \ldots \partial^{\beta_d}f_d)(x) ~\langle
\varphi ,\partial^\gamma \delta _{f(x)} \rangle. \eea For $\omega
\not\in \tilde{N}$, define $Y_t (\omega ):{\cal E}'\rightarrow {\cal
E}'$ by \bea Y_t(\omega )(\psi)&=& \sum\limits_{|\alpha |\leq
N}(-1)^{|\alpha |} \sum\limits_{\stackrel{|\gamma |\leq|\alpha
|}{|\beta _i |\leq |\alpha |}}
(-1)^{|\gamma |} \nonumber \\
&=& \int\limits_V g_\alpha (x)~ P_\gamma (\partial^{\beta _1} X_1
\ldots \partial^{\beta _d}X_d) (t,x,\omega )~  \partial^\gamma
\delta_{X(t,x,\omega )}~ dx~~ \eea Take $Y_t (\omega )=0$ if $\omega
\in \tilde{N}$.

Let ${\cal S}$ be the space of smooth rapidly decreasing functions
on $\R^d$  with dual ${\cal S}'$,  the space of tempered
distributions. It is well known (\cite{I}), that ${\cal S}$ is a
nuclear space, and that ${\cal S}=\displaystyle{\bigcap_{p>0} \left(
S_p,\| \cdot\|_p\right)}$, where the Hilbert spaces $S_p$ are
equipped with increasing norms $\| \cdot\|_p$, defined by the inner
products \Bea \langle f,g\rangle_p=\sum_{|k|=0}^\infty \left(
2|k|+d\right)^{2p}\left\langle f,h_k\right\rangle\left\langle
g,h_k\right\rangle,\,\,\,\,f,g\in {\cal S} \Eea
 Above, $\{h_k\}_{|k|=0}^\infty$ is an orthonormal basis for 
$L^2\left( \R^d,dx\right)$
given by Hermite functions (for $d=1$, $h_k(t) = \bigl(
2^kk!\sqrt{\pi}\bigr)^{-1/2}\exp\{-t^2/2\}H_k(t)$, with $H_k(t)$, a
Hermite polynomial. see \cite{I}), and $\langle\cdot ,\cdot \rangle$
is the usual inner product in $L^2\left( \ R^d,dx\right)$. We also
have ${\displaystyle{\mathcal S}'= \bigcup_{p>0}\left( S_{-p},\|
\cdot \|_{-p}\right)}$. Note that for all $-\infty<p<\infty$, $S_p$
is the completion of ${\mathcal S}$ in $\|\cdot \|_p$ and $S'_p$ is
isometrically isomorphic with $S_{-p}$, $p>0$. We then have the
following proposition.

\begin{Proposition} Let $\psi$ be a distribution with compact support having
representation (3.1). Let $p >0$ be such that $\partial^\alpha
\delta_x \in S_{-p}$ for $|\alpha | \leq N$. Then $(Y_t (\psi ))_{t
\geq 0}$ is an $S_{-p}$ valued continuous adapted process such that
for all $t \geq 0$,
$$ Y_t (\psi)=X_t^\ast (\psi) ~a.s.~ P.$$
\end{Proposition}

{\bf Proof:} From Theorem 4.6.5 of \cite{K}, it follows that for
any multi index $\gamma$ and compact set $K \subseteq \R^d$, $
\sup\limits_{s \leq T} ~\sup\limits_{x \in K} |\partial^\gamma X_i
(s,x,\omega )| < \infty$ a.s. for all $ T > 0$.

From this result, the fact that $\partial_{i} : S_{-{p-\frac{1}{2}}}
\rightarrow S_{-p}$ is bounded, and Theorem 2.1 of \cite{RT}, it
follows that for all $ T >0 $ a.s. \bea
\!\!\!\!\!\!\sup\limits_{t \leq T} \int\limits_V |g_\alpha (x)|
~|P_\gamma (\partial^{\beta _1} X_1, \ldots \partial^{\beta _d} X_d
) (t,x,\omega ))|~\|\partial^\gamma \delta_{X(t,x,\omega )}\|_{-p}
dx < \infty. \eea It follows that $(Y_t (\psi))$ is a well defined
$S_{-p}$-valued process. Since for all $x \in \R^d, (\partial^\gamma
X_i (t,x,\omega ))$ is an adapted process, and jointly measurable in
$(t,x,\omega )$ it follows that $(Y_t (\psi))$ is an $S_{-p}$-valued
adapted process.

To show that the map $ t \rightarrow Y_t(\psi) $ is almost surely
continuous in $S_{-p}$, we first observe that for any $p \in \R$
 and any $\phi \in S_p$, the map $ x \rightarrow \tau_x\phi : \R^d
 \rightarrow S_p $ is continuous, where $\tau_x : S' \rightarrow S'$
 is translation by $x \in \R^d $. To see this, let $x_n \rightarrow x \in \R^d $.
 Note that from Theorem 2 of \cite{RT}, given $ \epsilon > 0 $,
 there exists $ \psi \in S $ such that \Bea \| \tau_x\phi -
 \tau_{x_n}\phi \|_p &<&  \| \tau_x\psi -
 \tau_{x_n}\psi \|_p + \frac{\epsilon}{2}. \Eea
 From the definition of $\|.\|_p$ we have \Bea \| \tau_x\psi -
 \tau_{x_n}\psi \|_p^2 &=& \sum_{k} (2|k| + d)^{2p}\langle\tau_x\psi -
 \tau_{x_n}\psi, h_k \rangle^2 .\Eea Since $\psi \in S_q$ for every $q$ and since
 the  result is true for $p = 0$, the right hand side above tends to zero by
 dominated convergence theorem and \Bea \| \tau_x\phi -
 \tau_{x_n}\phi \|_p &<&   \epsilon \Eea for large $n$ thus proving
 the continuity of the map $ x \rightarrow \tau_x\phi $. In particular
 if $p$ and $\alpha$ are as in the statement of the theorem, the map
 $x \rightarrow \partial^\gamma \delta_x $ is continuous in $S_{-p}$ for
 $|\gamma| \leq |\alpha| $ .  Now the continuity of $ t \rightarrow Y_t(\psi)(\omega) $
 follows from the continuity in the $t$ variable of the processes
  $(\partial^\gamma X_i (t,x,\omega ))$ and $(\partial^\gamma
  \delta_{X(t,x,\omega)})$
and the dominated convergence theorem .

We now verify that $Y_t (\psi) =X_t^\ast (\psi)$. Let $\varphi \in
C^\infty (\R^d)$. From (3.3), \Bea \langle Y_t (\omega ) \psi, \phi
\rangle &=& \sum\limits_{|\alpha | \leq N} (-1)^{|\alpha |}
\sum\limits_{\stackrel{|\gamma | \leq |\alpha |}{|\beta _i |
\leq |\alpha |}} (-1)^{|\gamma |}\\
&& \int\limits_V g_\alpha (x) P_\gamma
(\partial^{\beta _1} X_1, \ldots \partial X_d^{\beta _d} ) (t,x,\omega )
~\langle \partial^\gamma \delta_{X(t,x,\omega )}, \phi  \rangle~ dx\\
&=& \sum\limits_{|\alpha | \leq N} (-1)^{|\alpha |} \int\limits_V
g_\alpha (x) \partial^\alpha (\phi \circ X (t,x,\omega )) dx ~~~~~\mbox{(by (3.2))}\\
&=& \sum\limits_{|\alpha | \leq N} \langle \partial^\alpha g_\alpha ,\phi \circ
X(t,\cdot,\omega ) \rangle ~~~~~~~\mbox{(by (3.1))}\\
&=& \sum\limits_{|\alpha | \leq N} \langle \partial^\alpha g_\alpha
, X_t (\omega )
\phi \rangle \\
&=& \langle \psi , X_t (\omega ) \phi \rangle.
\Eea
\hfill${\Box}$

Recall that $C^\infty$ denotes the space of smooth functions on
$\R^d$. Define the operators $A : C^\infty \rightarrow
L(\R^r,C^\infty)$ and $L: C^\infty \rightarrow C^\infty$ as follows:
For $\varphi \in C^\infty $, $  x \in \R^d$, \Bea A\varphi &=&
(A_1\varphi,\cdots,A_r\varphi),\\A_i\varphi (x) &=&
\sum\limits^d_{k=1} \sigma_{i}^k(x)
\partial_k \varphi(x),
\\L\varphi (x) &=& \frac{1}{2} \sum\limits^d_{i,j=1} (\sigma \sigma
^t)_{j}^i (x) ~\partial^2_{ij} \varphi (x)  + \sum\limits^d_{i=1}
b^i (x) ~\partial _i \varphi (x) . \Eea We define the adjoint
operators $A^\ast : {\cal E}' \rightarrow L(\R^r, {\cal E}')$ and
$L^\ast : {\cal E}' \rightarrow {\cal E}'$ as follows: \Bea
A^\ast\varphi &=&
(A^\ast_1\varphi,\cdots,A^\ast_r\varphi),\\
A^\ast_i \psi &=& - \sum\limits^d_{k=1}
\partial_k
(\sigma_{i}^k \psi), \\
L^\ast \psi &=& \frac{1}{2} \sum\limits^d_{i,j=1} \partial^2_{ij}
((\sigma \sigma^t)_{j}^i \psi)  - \sum\limits^d_{i=1} \partial_i
(b^i \psi). \Eea The following proposition gives the boundedness
properties of $L^\ast$ and $A^\ast$. For $K \subset \R^d$, let ${\cal
E}' (K) \subseteq {\cal E}'$ be the subspace of distributions whose
support is contained in $K$.

\begin{Proposition} Let $p >0 $ and $ q > [p]+2$, where $[p]$ denotes the
largest integer less than or equal to $p$. Then,  $A^\ast : S_{-p}
\cap {\cal E}' (K) \rightarrow L(\R^r, S_{-q} \cap {\cal E}' (K))$
and $L^\ast : S_{-p} \cap {\cal E}' (K) \rightarrow S_{-q} \cap
{\cal E}' (K)$. Moreover, there exists constants $C_1 (p) >0, C_2
(p) >0$ such that \Bea
\|A^\ast \psi\|_{HS(-q)} \leq  C_1 (q) ~\|\psi\|_{-p}, ~~~~
\|L^\ast \psi\|_{-q} \leq C_2 (q) ~\|\psi\|_{-p} \Eea
where \Bea 
 |A^\ast \psi\|^2_{HS(-q)}=  \sum\limits^r_{i=1}
\left\| \sum\limits^d_{k=1}
\partial_k (\sigma_{i}^k \psi ) \right\|^2_{-q} = 
 \sum\limits^r_{i=1}\|A_i^\ast\psi\|_{-q}^2. \Eea
\end{Proposition}

{\bf Proof} Clearly $A^\ast : {\cal E}' (K) \rightarrow L(\R^r,
{\cal E}' (K))$ and $L^\ast : {\cal E}' (K) \rightarrow {\cal
E}'(K)$. We prove the bounds for $A^\ast$. The bounds for $L^\ast$
follow in a similar fashion. By definition, if $q \geq p+1/2$, \bea
\|A^\ast \psi \|^2_{HS(-q)} &=& \sum\limits^d_{i=1}
~\|\sum\limits^d_{k=1}~
\partial_k (\sigma _{ki} \psi )\|^2_{-q} \nonumber \\
&\leq& C ~\sum\limits^d_{i=1} ~\sum\limits^d_{k=1}~
\|\sigma_{ki}\psi \|^2_{-(q-\frac{1}{2})}. \eea Let $\sigma$ denote a
$C^\infty$ function. We first show that the map $\psi \rightarrow
\sigma \psi :{\cal S}_n \cap {\cal E}' (K)\rightarrow {\cal S}_n
\cap {\cal E}' (K)$ satisfies \Bea \|\sigma \psi \|_n \leq C~ \|\psi
\|_n \Eea where the constant $C$ depends on $\sigma $, $K$ and $n$.
First assume that $\psi \in {\cal S} \cap {\cal E}'(K)$. We then
have (see \cite{RT}, Proposition 3.3b ) \Bea \|\sigma \psi \|^2_n
\leq C_1 \sum\limits_{|\alpha |+|\beta |\leq 2n} \|x^\alpha
\partial^\beta (\sigma \psi )\|^2_0. \Eea Clearly, \Bea \|x^\alpha
\partial^\beta (\sigma \psi )\|^2_0 \leq C_2 \sum\limits_ {|\gamma|
\leq |\beta |} \|x^\alpha \partial^\gamma \psi \|^2_0. \Eea It
follows that \Bea \|\sigma \psi \|^2_n &\leq& C_3
\sum\limits_{|\alpha |+|\gamma | \leq
2n} \|x^\alpha \partial ^\gamma \psi \|^2_0 \\
&\leq& C  \|\psi \|^2_n. \Eea We now extend this to $\psi \in
{\cal S}_n \cap {\cal E}' (K)$: Since ${\cal S}$ is dense in ${\cal
S}_{n}$, we can get $\psi _m \in {\cal S}, \psi _m \rightarrow \psi$
in ${\cal S}_n$. By multiplying by an appropriate $C^\infty$
function with compact support can assume, $\psi _m \in {\cal S} \cap
{\cal E}' (K)$. By the above inequality applied to $\psi_m$, it
follows that $\sigma \psi _m$ converges in ${\cal S}_n$, and hence
weakly to a limit $\varphi$ in ${\cal S}'$. Hence if $f$ is
$C^\infty$ with compact support, \Bea \langle \varphi ,f\rangle &=&
\lim\limits_{m \rightarrow \infty}
\langle \sigma \psi _m, f \rangle \\
&=& \lim\limits_{m \rightarrow \infty} \langle \psi _m, \sigma f\rangle\\
&=& \langle \psi , \sigma f \rangle = \langle \sigma \psi , f
\rangle. \Eea Hence $\sigma \psi _m \rightarrow \sigma \psi$ in
${\cal S}_n$ and the above inequality  follows for $\psi \in {\cal
S}_n \cap {\cal E}'(K)$.

Now suppose $\psi \in {\cal S} \cap {\cal E}' (K)$. \Bea \|\sigma
\psi \|_{-n} &=& \sup\limits_{\stackrel{\|\varphi\|_n \leq 1}
{\varphi \in {\cal S}}}
| \langle \sigma \psi , \varphi \rangle | \\
&=& \sup\limits_{\stackrel{\|\varphi \|_n \leq 1} {\varphi \in {\cal
S}}} ~|\langle g \sigma \psi , \varphi \rangle | \Eea where $g \in
C^\infty, g =1$ on $K$, $\mbox{supp} (g) \subseteq K^\epsilon$, an
$\epsilon$-neighbourhood of $K$. Therefore, \bea \|\sigma \psi
\|_{-n} &\leq& \|\psi \|_{-n} ~\sup\limits_{\stackrel{\|\varphi \|_n
\leq 1}{\varphi \in {\cal
S}}} \|\sigma g\varphi \|_n \nonumber\\
&\leq& C  \|\psi \|_{-n}. \eea In the same way as for $n \geq
0$, we can extend the above inequality to $\psi \in {\cal S}_{-n}
\cap {\cal E}'(K)$. Now the proof can be completed using (3.5) and
(3.6) and by choosing a $q \in \R $ such that $ q \geq m +
\frac{1}{2}
> m \geq p $ for some integer $m$. In particular, we may take $q >
[p] + 2$. \hfill${\Box}$

\begin{Theorem} Let $\psi \in {\cal E}'$ have the representation (3.1).
Let $p >0$ be such that $\partial^\gamma \delta_x \in S_{-p},
|\gamma | \leq N$. Let $q >p$ be as in Proposition 3.2. Then the
$S_{-p}$-valued continuous, adapted process $(Y_t (\psi))_{t \geq
0}$ satisfies the following equation in $S_{-q}$: a.s., \bea Y_t
(\psi) =\psi +\int\limits^{~~t}_{0} A^\ast (Y_s (\psi)) \cdot dB_s
+\int\limits_0^{~~t} L^\ast (Y_s (\psi))~ds \eea for all $t \geq 0$.
\end{Theorem}

{\bf Proof:} From Proposition 3.2 and the estimate (3.4), it follows that for
$t \geq 0$, a.s.
$$ \int\limits_0^{~~t} \|A^\ast Y_s (\psi)\|^2_{HS(-q)} ~ds+\int\limits_0^{~~t}
\|L^\ast \psi_s (\psi) \|_{-q} ds < \infty.$$ Hence the right hand side 
of (3.7)
is well defined. Let $\varphi \in C^\infty$. Then by Ito's formula
\Bea X_t( \varphi) &=& \varphi +\int\limits_0^{~~t} X_s (A\varphi )
\cdot dB _s   + \int\limits_0^{~~t} X_s (L \varphi ) ds \Eea where
the  integrals on the right hand side are understood as $C^\infty (\R^d)$ 
valued
processes given by $(t,x,\omega ) \rightarrow \int\limits_0^{~~t}
A\varphi (X(s,x,\omega ) \cdot dB_s$ and $(t,s,\omega ) \rightarrow
\int\limits_0^{~~t} L\varphi (X(s,x,\omega )) ds$. Since $\psi$ has
compact support, by localising, we may assume that these processes
have their supports contained in the support of $\psi, $ a fixed
compact set not depending on $t$ and $\omega$. In particular, they
belong to $S \subset S_p $. Using Proposition 3.1, we get, for
$\varphi \in S$, \Bea
\langle Y_t (\psi), \varphi \rangle &=& \langle \psi, X_t (\varphi) \rangle \\
&=& \langle \psi, \varphi \rangle + \int\limits_0^{~~t} \langle \psi,
X_s (A\varphi ) \rangle \cdot dB_s  + \int\limits_0^{~~t} \langle \psi, X_s (L\varphi ) \rangle~ ds \\
&=& \langle \psi, \varphi \rangle +\int\limits_0^{~~t} \langle A^\ast
Y_s (\psi), \varphi \rangle \cdot dB_s + \int\limits_0^{~~t} \langle L^\ast Y_s (\psi), \varphi \rangle ~ds\\
&=& \langle \psi +\int\limits_0^{~~t} A^\ast Y_s (\psi) \cdot dB_s +
\int\limits^{~~t}_0 L^\ast Y_s (\psi)~ ds, ~~\varphi \rangle \Eea
and the result follows. In the above calculations we have used the
fact that $T\int\limits_0^t A^\ast Y_s (\psi) \cdot dB_s=
\int\limits_0^t TA^\ast (Y_s (\psi)) \cdot dB_s$ for any bounded
linear functional $T: S_{-q} \rightarrow \RR$. \hfill${\Box}$

\section{{\protect\large {\bf Probabilistic representations}}}
\setcounter{equation}{0}

In this section we prove the probabilistic representations of
solutions to the initial value problem for the parabolic operator
$\partial_t - L^\ast$. We also show uniqueness for the solutions
of the initial value problem under the `Monotonicity conditions'.
We first prove some estimates on $\|\delta_x\|_{-p}$ and
$\|\partial^\gamma \delta_x \|_{-p}$ that are required later.

\begin{Theorem} $\delta_x \in S_{-p}$ iff $ p > \frac{d}{4}$.
Further if $p > \frac{d}{4}$, then
$$\lim\limits_{|x|\rightarrow \infty}\|\delta_x\|_{-p}=0.$$
 In particular, if for a multi-index $\gamma \in \ZZ_+^d$,
$p > \frac{ d}{4} +\frac{|\gamma |}{2}$, then
$$\sup\limits_{x \in \R^d} \| \partial^\gamma \delta_x\|_{-p} < \infty.$$
\end{Theorem}

{\bf Proof:}  Since $ \partial^\gamma : S_{-q} \rightarrow S_{-q-\frac
{|\gamma|}{2}}$ is a continuous operator we have for $p > \frac{d}{4} +
\frac{|\gamma|}{2}$,
$$  \| \partial^\gamma \delta_x\|_{-p} \leq C \|\delta_x\|_{(-p -\frac
{|\gamma|}{2})}. $$ Further from Theorem (2.1) of \cite{RT}, it
follows that for any compact set $ K $ contained in $\R^d$,
$\sup\limits_{x \in K} \|
 \delta_x\|_{-(p+ \frac{|\gamma|}{2} )} < \infty $.  Since
$p+ \frac{|\gamma|}{2} > \frac{d}{4}$, the last statement of the theorem
now follows from the second statement.

The proof of the first part of the theorem uses the generating
function for Hermite functions given by Mehler's formula (see
\cite{T}, page 2). First we note that
$$(2n+d)^{-2p} =\frac{1}{(2p-1)!} \int\limits^{~~\infty}_0 t^{2p-1}~
e^{-(2n+d)t} dt.$$
Hence,
\Bea
\|\delta_x \|^2_{-p} &=& \sum\limits^\infty_{n=0} (2n+d)^{-2p} ~
\sum\limits_{|k|=n} |h_k (x)|^2 \\
&=& \frac{1}{(2p-1)!} \int\limits^{~~\infty}_0 t^{2p-1} ~g(t,x) dt \Eea        where $$    g(t,x)  \sum\limits^\infty_{n=0} e^{-(2n+d)t}
\sum\limits_{|k|=n} |h_k (x)|^2. $$   Using Mehler's formula , \Bea
g(t,x) &=& e^{-dt} \pi^{-\frac{d}{2}} (1-e^{-4t})^{-\frac{d}{2}} \\
&& \times ~e^{-\frac{1}{2}(\frac{1+e^{-4t}}{1-e^{-4t}})2|x|^2 +
(\frac{e^{-2t}}{1-e^{-4t}})2|x|^2}\\
&=& e^{-dt} \pi^{-\frac{d}{2}} (1-e^{-4t})^{-\frac{d}{2}} \\
&& \times ~e^{- (tanh~ t)|x|^2}. \Eea It is easy to see that for all $ x $
\Bea
&& g(t,x) \sim (1-e^{-4t})^{-\frac{d}{2}}, ~~t \rightarrow 0 \\
\mbox{and}~ && g(t,x) \sim e^{-dt} e^{-(tanh~t)|x|^2}, ~~t
\rightarrow \infty \Eea where $a(t) \sim b(t) \leftrightarrow
\frac{a(t)}{b(t)}\rightarrow c~(>0)$. It follows that for all $ \epsilon > 0 $,
$$\int\limits_0^{~~\epsilon} t^{2p-1} g(t,x) dt < \infty $$
iff $p >\frac{d}{4}$. Also $g(t,x)\rightarrow 0$ as $|x|\rightarrow
\infty$, for every $t >0$. Hence by the dominated convergence theorem,
\Bea
\lim\limits_{|x|\rightarrow \infty} \|\delta_x \|^2_{-p} &=& \lim
\limits_{|x|\rightarrow \infty} \frac{1}{(2p-1)!} \int\limits_0^{~~\infty}
t^{2p-1} g(t,x) dt \\
&=& 0.
\Eea
\hfill$\Box$

\begin{Proposition} Let $\psi \in {\cal E}'$ with
representation (3.1). Let $ p > \frac{d}{4} + \frac{N}{2}$ where $N
= order (\psi) + 2d$. Let $(Y_t(\psi))$ be the $S_{-p}$ valued
continuous adapted process defined by (3.3). Then for all $ T >0$,
$$\sup\limits_{t \leq T} E\|Y_t (\psi)\|^2_{-p} < \infty.$$
\end{Proposition}

{\bf Proof}: Using the representation (3.3) for $Y_t(\psi)$ and the
result of Theorem 4.1,
$$E\|Y_t (\psi)\|^2_{-p} \leq C \cdot \sum\limits_{|\alpha |\leq N} \sum\limits_
{\stackrel{|\gamma | \leq |\alpha |}{|\beta _i| \leq |\alpha
|}}\int_V E \{P_\gamma (\partial^{\beta _1} X_1, \ldots
\partial^{\beta _d} X_d) (t,x,\omega ) \}^2 dx.$$ Hence it suffices to
show,
$$\sup\limits_{t \leq T} \sup\limits_{x \in V} E\{P_\gamma (\partial^{\beta _1}
X_1, \ldots \partial^{\beta _d} X_d) (t,x,\omega ) \}^2 < \infty.$$
This follows from Theorem 2.2 and completes the proof.
\hfill$\Box$

We now consider solutions to the following initial value problem
\bea
\left. \begin{array}{lll} \frac{\partial \psi_t}{\partial t} &=&
L^\ast \psi_t \\
\psi_0 &=& \psi \end{array} \right\}
\eea
for $\psi \in {\cal E}'$. By a solution to (4.1) we mean an
 $S_{-p}$-valued
continuous function $\psi:[0,T]\rightarrow S_{-p}$ for some
$p >0$, such that the equation
$$\psi_t =\psi +\int\limits^{~~t}_0 L^\ast \psi_s ~ds$$
holds in $S_{-q}$, where $q > p$ is such that $L^\ast \psi_s \in S_{-q},
0 \leq s \leq T $ and is Bochner integrable on $[0,T]$ with
respect to Lebesgue measure. Note that, in our definition the
initial value $\psi$ belongs to $S_{-p}$. Note also that we can take
$q >p$ given by Proposition 3.2 .

Define $\psi_t = EY_t (\psi), t \geq 0$. If $\psi$ has
compact support, then the previous proposition implies
that $\psi_t$ is a well defined element of an appropriate
Hilbert space $S_{-p}$. The following theorem gives us
the existence of solutions of equation (4.1) and a
stochastic representation of its solutions.

\begin{Theorem}
Let $\psi \in {\cal E}' $ have representation (3.1) and let $p$ be as in
Proposition 4.2. Then $(\psi_t)_{0 \leq t \leq T}$ is an
$S_{-p}$-valued solution of (4.1).
\end{Theorem}

{\bf Proof }: Let $ q > p$ be as in Proposition 3.2. Because of
Proposition (4.2) we can take expectations in (3.7): \Bea
E\|\int\limits_0^{~~t} A^\ast Y_s (\psi) \cdot dB_s\|^2_{-q} &\leq&
E\int\limits_0^{~~t} \|A^\ast Y_s (\psi)\|^2_{HS(-q)} ds \\
&\leq& C \cdot E \int\limits_0^{~~t} \|Y_s (\psi )\|^2_{-p} ds \\
&<& \infty .
\Eea
Similarly,
\Bea
E\left\| \int\limits_0^{~~t} L^\ast Y_s (\psi) ds \right\|_{-q} &\leq&
E \int\limits^{~~t}_0 \|L^\ast Y_s (\psi)\|_{-q} ds \\
&\leq& C  E \int\limits_0^{~~t} \|Y_s (\psi)\|_{-p} ds <
\infty. \Eea In particular, $E\int\limits_0^{~~t} A^\ast Y_s (\psi)
\cdot dW_s =0$. Taking expectations in (3.7), we get
$$\psi_t =S_t (\psi) =EY_t (\psi) =\psi +\int\limits^{~~t}_0 EL^\ast
Y_s(\psi) ds.$$ Since $L^\ast : S_{-p} \cap {\cal E}'(K)
\rightarrow S_{-q} \cap {\cal E}' (K)$ is a bounded operator
(Proposition 3.2), $EL^\ast Y_s (\psi)=L^\ast EY_s (\psi) =
L^\ast \psi_s.$
\hfill$\Box$

We now consider the uniqueness of solutions to the initial value
problem (4.1). For $p >0$, let $q >p$ be as in Proposition 3.2 ($q >
[p] + 2 $) so that $A^\ast : S_{-p} \cap {\cal E}' \rightarrow
L(\R^r, S_{-q} \cap {\cal E}')$ and $L^\ast : S_{-p} \cap {\cal E}'
\rightarrow S_{-q} \cap {\cal E}'$ are bounded operators. The pair
$(A^\ast, L^\ast)$ is said to satisfy the monotonicity inequality in
$S_{-q} \cap {\cal E}'$ if and only if there exists a constant
$C=C(p)$, such that, \bea 2 \langle \varphi,L^\ast \varphi
\rangle_{-q} +\|A^\ast \varphi\|^2_{HS(-q)} \leq C
~\|\varphi\|^2_{-q} \eea holds for all $\varphi \in S_{-p} \cap
{\cal E}'$.

\begin{Theorem} Let $\psi \in {\cal E}'(K)$ have representation 3.1. Let
 $p > \frac{d}{4} +\frac{N}{2}$ and $N=$ order $\psi +2d$ (In particular
 $\psi \in S_{-p} )$. Let $ q > [p] + 2 $.  Suppose the
pair $(A^\ast, L^\ast)$ satisfies (4.2). Then, the initial value
problem (4.1) has a unique $S_{-p}$ valued solution given by
$\psi_t=EY_t (\psi)$.
\end{Theorem}

{\bf Proof:} The existence has been proved. It suffices to show
uniqueness. Let $\psi_t'$ be another $S_{-p}$-valued solution. Let
$\varphi_t =\psi_t-\psi_t'$. Then $(\varphi_t)$ satisfies, \Bea
\varphi_t =\int\limits_0^{~~t} L^\ast \varphi_s ~ds, ~0 \leq t
\leq T \Eea in $S_{-q}$ for $q > [p]+2$. Hence \Bea
\|\varphi_t\|^2_{-q} &=& 2 \int\limits_0^{~~t} \langle \varphi_s,
L^\ast \varphi_s \rangle_{-q} ~ds \\
&\leq& \int\limits_0^{~~t} \{2 \langle \varphi_s ,L^\ast \varphi_s
\rangle_{-q} +\|A^\ast \varphi_s\|^2_{HS (-q)} \} ds \\
&\leq& C \cdot \int\limits^{~~t}_0 \|\varphi_s\|^2_{-q} ~ds.\Eea
Now, the Gronwall inequality implies $\varphi_t \equiv 0, 0 \leq t
\leq T$. $\hfill{\Box}$

{\bf Remark:} When $\sigma_{j}^i (x)$ and $b^i(x)$ are constants
(independent of $x$), the inequality (4.2) for $(A^\ast,L^\ast)$ was
proved in \cite{GMRMon}.

We now show that the expected value of the `stochastic fundamental
solution' $(Y_t(\delta_x))=(\delta_{X(t,x)})$ of equation (3.5)
gives us the fundamental solution of (4.1). Let $P(t,x,A)=P(X(t,x)
\in A)$ be the transition function of the diffusion $(X(t,x))$.

\begin{Theorem} Let $\psi$ be an integrable function with compact
support . Let $ p > \frac{d}{4}+1$.

a)  Then for $0 \leq t \leq T$ and $x \in \R^d,  P(t,x,\cdot) \in
S_{-p}$ and \Bea P(t,x,\cdot) = E\delta_{X(t,x)} = EY_t(\delta_x)
\Eea

b)Let $(\psi_t)$ be the $S_{-p}$-valued solution of equation (4.1)
given by $\psi_t =EY_t(\psi)$. Then \Bea \psi_t = \int \psi(x)
~P(t,x,\cdot)dx \Eea where the  integral in the right hand side
above is an $S_{-p}$ valued Bochner integral.
\end{Theorem}

{\bf Proof:} a) The equality $E\delta_{X(t,x)}=EY_t(\delta_x)$
follows from the fact that if $\phi \in S$, \Bea \langle EY_t
(\delta_x), \phi \rangle &=& E
\langle Y_t (\delta_x), \phi \rangle  \\
&=& E \langle \delta_x, \phi (X(t,\cdot)) \rangle \\
&=& E \phi (X(t,x)) =E \langle \delta_{X(t,x)}, \phi \rangle \\
&=& \langle E\delta_{X(t,x)}, \phi \rangle. \Eea On the other hand,
we have $E \phi(X(t,x))=\langle P (t,x,\cdot), \phi \rangle$ and the
equality asserted in a) follows. To show that $ P(t,x,\cdot)$
belongs to $S_{-p}$, we first note that using Theorem (2.1) of
\cite{RT}, there exists a polynomial $P(x)$ such that \Bea \|EY_t
(\delta_x)\|_{-p} &=& \|E
\delta_{X(t,x)}\|_{-p} \\
&\leq& E\|\delta_{X(t,x)}\|_{-p} \\
&\leq& (E|P(X(t,x))|)~\|\delta_0\|_{-p}. \Eea It follows that
$P(t,x,\cdot)$ belongs to $S_{-p}$.

b) To prove b) we first show the Bochner integrability of
$P(t,x,\cdot)$ in $S_{-p}$. This follows from the calculations above
for part a)  and the fact that by Theorem 2.2, \Bea \sup\limits_{x
\in ~\mbox{supp} ~\psi} \|P(t,x,\cdot)\|_{-p} &\leq& \sup\limits_{x
\in ~\mbox{supp}
~\psi} (E|P(X(t,x))|)~\|\delta_0\|_{-p}\\
&<& \infty.
\Eea

Hence the integral $\int \psi(x)~P(t,x,\cdot)~dx$ is well defined
and belongs to $S_{-p}$. Further since $\psi$ is an integrable
function with compact support, the representation of $Y_t (\psi)$
given by equation (3.3) reduces to \Bea Y_t(\psi) =\int\limits_V
\psi(x) ~\delta_{X(t,x)} ~dx \Eea where supp$(\psi) \subset V$, $V$
an open set with compact closure. Hence \Bea \psi_t &=& E~Y_t (\psi)
=
\int\limits_V \psi(x) ~E\delta_{X(t,x)} ~dx \\
&=& \int\psi(x) ~P(t,x, \cdot)~dx. \Eea

$\hfill{\Box}$

Suppose now that $(X(t,x))$ has a density $p(t,x,y)$, i.e. \Bea
P(t,x,A) =\int\limits_A p(t,x,y) ~dy. \Eea We shall assume for the
rest of the section the following integrability condition on
$p(t,x,y)$ : For every compact set $K \subset \R^d$, \Bea
\int\!\!\int_{K \times K} p(t,x,y)~dx ~dy < \infty .\Eea

{\bf Corollary:} Let $K \subseteq \R^d$ be a compact set and $\psi$
be as in Theorem 4.5.
 If supp$(\psi) \subseteq K $, then
$\psi_t =EY_t (\psi)$ is given by a locally integrable function
$f(y)$ where $f(y)=\int \psi(x)~p(t,x,y)~dx$.

{\bf Proof:} Let $\phi \in C_c^\infty (\R^d)$. Then \Bea \langle f,
\phi \rangle &=& \int f(y) ~\phi (y)~dy = \int  ~\phi
(y)~(\int \psi (x) ~p(t,x,y)~dx)dy \\
&=& \int  ~\psi(x) (\int  ~p(t,x,y)~\phi (y)~dy)dx \\
&=& \int  ~\psi(x) ~\langle P(t,x,\cdot), \phi \rangle dx \\
&=& \left\langle \int \psi(x)~P(t,x,\cdot)~dx , \phi \right\rangle \\
&=& \langle \psi_t, \phi \rangle \Eea where the third equality
follows from our assumptions on $p(t,x,y)$ and Fubini's theorem.
$\hfill{\Box}$

We now consider the self adjoint case $L^\ast = L $. We deduce,
under some mild integrability conditions, the well known result
that the transition density is symmetric.

\begin{Theorem} Suppose $\sigma_j^i, b^i$ are $C^\infty$,
bounded with bounded derivatives of all orders. Suppose $(A^\ast,
L^\ast)$ satisfy the monotonicity condition (4.2). Suppose further
that $L^\ast=L$ . Then for $0 < t \leq T, p(t,x,y)=p(t,y,x)$ for
every $(x,y)$ outside a set of zero Lebesgue measure in $\R^d \times
\R^d$.
\end{Theorem}

{\bf Proof:} Let $\psi \in C_c^\infty (\R^d)$. Let $\psi (t,x)
=E\psi (X(t,x))$. Under our assumptions on $\sigma_j^i, b^i$ and
$\psi$, and the assumption $L = L^\ast$, it is well known (see
\cite{B}, p.47 ) that $\psi (t,x)$ is a classical solution of the
initial value problem (4.1). On the other hand by the monotonicity
condition (4.2), we have
 uniqueness of  the initial value problem
(Theorem 4.4) and hence by Corollary to Theorem 4.5, we have \Bea
\int \psi (y)~p(t,x,y) ~dy &=&
E\psi (X_t^x) = \psi (t,x) \\
&=& \int \psi (y)~p(t,y,x)~dy \Eea for a.e. $x \in \R^d$. Since
$\psi$ is arbitrary, this implies $p(t,x,y)=p(t,y,x)$ for every
$(x,y)$ outside a set of zero Lebesgue measure and completes the
proof. $\hfill{\Box}$

In the constant coefficient case we can deduce the following well
known result.

\begin{Proposition} Suppose $\sigma_j^i, b^i$ are constants and
that $(A^\ast, L^\ast)=(A,L)$ satisfy the monotonicity inequality
(4.2). Then $p(t,x,y)=p(t,0,y-x)$ for almost every $(x,y)$ with
respect to the Lebesgue measure on $\R^d \times \R^d$.
\end{Proposition}

{\bf Proof:} Let $\psi$ be a continuous function with compact
support. By Corollary to Theorem 4.5, the distribution $EY_t (\psi)$
is given by the locally integrable function $\psi(t,y)=\int \psi (x)
~p(t,x,y)~dx$. On the other hand, by the uniqueness of solutions to
the SDE (3.7) (see \cite{GMR},\cite{GMRMon} ) we have a.s.
$Y_t(\psi)=\tau_{X_t} (\psi)$ for all $t \geq 0$, where $(X_t)$ is
the diffusion $(X(t,x))$ with $x=0$. In particular
$EY_t(\psi)=E\tau_{X_t}(\psi)$. The latter distribution is given by
a locally integrable function of $ y $.  We then have for a.e.
$y$, \Bea
\int \psi(x)~p(t,x,y)~dx &=& \psi (t,y) =EY_t (\psi) = E\tau_{X_t}(\psi)\\
&=& \int \tau_x \psi (y)~p(t,0,x)~dx \\ &=&\int \psi
(x)~p(t,0,y-x)~dx. \Eea Since $\psi $ is arbitrary, it follows that
$p(t,x,y)=p(t,0,y-x)$ for almost every $(x,y)$ with respect to
Lebesgue measure in $\R^d \times \R^d$. $\hfill{\Box}$

Define $S_t (\psi) =EY_t (\psi), t \geq 0$ for $\psi \in {\cal E}'$.
Then $ S_t : {\cal E}' \rightarrow S' $. Let $(T_t)_{t \geq 0}$ be
the semigroup corresponding to the diffusion $(X_t)$ solving (2.1)
i.e. for $f \in {\cal S}, T_t f(x)=Ef(X(t,x ))$. We have the
following result:

\begin{Theorem} a) $T_t : {\cal S} \rightarrow C^\infty$ and we have $S_t
=T_t^{\ast}$ in the sense that
$$\langle S_t (\psi), \phi \rangle = \langle \psi, T_t \phi\rangle$$
for all $\psi \in {\cal E}'$ and $\phi \in {\cal S}$.

b) Let $K \subset \R^d$ be a compact set and $p >0$. Then for $ q >
\frac{5}{4}d+ [p]+1$, $S_t:S_{-p} \cap {\cal E}' (K) \rightarrow
S_{-q}$ is a bounded linear operator. Further, for any $T >0$,
$$ \sup\limits_{t \leq T} \|S_t \|_H < C(T)$$
where $\| \cdot \|_H$ is the operator norm on the Banach space
$H$ of bounded linear operators from $S_{-p} \cap {\cal E}'(K)$ to
$S_{-q}$.
\end{Theorem}

{\bf Proof:} a) Clearly for $f \in S$, that $T_tf(x) = Ef(X(t,x))$
is ~$\C^\infty$ follows from Theorem 2.2 and the dominated
convergence theorem. In other words, $T_t : S \rightarrow C^\infty
$. Also if $\psi \in {\cal E}'$ and $ p > \frac{d}{4}+\frac{N}{2}, N
= order(\psi) +2d$, then Proposition 4.2 implies $EY_t(\psi) =
S_t(\psi) \in S_{-p} \subset S'$ . Hence $ S_t : {\cal E}'
\rightarrow S'$. We then have for $\phi \in S, \psi \in {\cal E}'$,
\Bea \langle S_t (\psi), \phi \rangle &=& E \langle Y_t (\psi), \phi
\rangle = E \langle \psi, X_t (\phi)\rangle \\
&=& E \langle \psi, \phi'(\phi(X(t,\cdot))) \rangle = \langle \psi,
\phi'E(\phi(X(t,\cdot))) \rangle \\
&=& \langle \psi, T_t(\phi) \rangle \Eea where $\phi' \in
C_c^\infty, \phi'=1$ on support of  $\psi$. Here the last but one
equality follows from the fact that $ E\phi'(\phi(X(t,\cdot)))
=\phi'E\phi(X(t,\cdot))$, and the fact proved below that $E
\|\phi'\phi(X(t,\cdot))\|_p \leq E\|\phi'(\phi(X(t,\cdot))\|_{[p]+1}
< \infty$ for every $p> 0$ and hence $ E\phi'(\phi(X(t,\cdot))$
belongs to $S_p$ for every $p
>0$ and in particular belongs to $S$. This completes the proof of
part a).

b) Without loss of generality, we may assume $p$ is an integer
(since $S_{-p} \subseteq S_{-([p]+1)}$, where $[p]=$ greatest
integer $\leq p$). Note that $\psi \in S_{-p} \cap {\cal E}'(K)$
implies order $\psi \leq 2p$ - this follows from the fact that if
$p$ is an integer and if support of $\phi$ is contained in $K$ then
(see \cite{RT}) \Bea \| \phi\|^2_p \leq C_1 \sum\limits_{|\alpha |
+|\beta | \leq 2p} \|x^\alpha \partial^\beta (\phi)\|^2_0 \leq C_2\|
\phi \|^2_p. \Eea Hence if $q$ is as in the statement of the theorem
and N = order of $\psi + 2d$, then $q > \frac{d}{4}+\frac{N}{2}$,
which in turn
 implies, by
Proposition 4.2, that $E\|Y_t (\psi)\|^2_{-q} < \infty$.
 In particular, $S_t (\psi) \in S_{-q}$ and hence it suffices to show
that there exists  $ C=C(T,p,K) >0$ such that for $\psi \in S_{-p} \cap
{\cal E}' (K), \phi \in S$, \Bea |\langle S_t (\psi), \phi \rangle |
\leq C\|\psi\|_{-p} \|\phi\|_p .\Eea Note that the left hand side
above is given by \Bea |\langle S_t (\psi), \phi \rangle | &=& |E
\langle Y_t (\psi), \phi
\rangle | = |E \langle \psi, X_t (\phi)\rangle |\\
&=& |E \langle \psi, \phi'X_t (\phi) \rangle | \Eea where $\phi' \in
C_c^\infty, \phi'=1$ on $K$.
Hence from the above equality we get
$$| \langle S_t (\psi) ,\phi \rangle | \leq \|\psi\|_{-p} (E\|\phi' X_t
(\phi)\|_p^2)^{\frac{1}{2}}.$$ Since (see \cite{RT})
$$\|\phi' X_t( \phi)\|^2_p \leq C \sum\limits_{|\alpha | +|\beta |
\leq 2p} \|x^\alpha \partial^\beta (\phi'X_t(\phi))\|^2_0$$
it suffices to show that for $|\alpha |+|\beta |\leq 2p$, there
exists a constant $C' > 0$ depending only on $T,p$ and $K$ such
that

\Bea E\|x^\alpha \partial^\beta (\phi'X_t(\phi))\|_0^2 &\leq& C'
\|\phi\|_p^2 .\Eea Let the support of $\phi'=K'$. We first compute the
expression inside the expectation sign in the right hand side above
for $\omega \notin \tilde{N} $ where $\tilde{N} $ is as in equation (3.3)
: \Bea \|x^\alpha \partial^\beta (\phi' X_t(\omega)(\phi))\|^2_0 &=&
\int\limits_{K'}
|x^\alpha \partial^\beta (\phi'(x) \varphi (X(t,x,\omega))|^2 dx \\
&=& \int\limits_{K'} |x^\alpha \sum\limits_{|\gamma |+|\gamma
'|=|\beta |}
\partial^\gamma \phi' (x) ~\partial^{\gamma '} \varphi (X(t,x,\omega)) |^2
dx \\
&=& \int\limits_{K'} | x^\alpha \sum\limits_{|\gamma |+|\gamma
'|=|\beta |}
\partial^\gamma \phi'(x) \sum\limits_{|\gamma'_1|,|\gamma'_2| \leq
|\gamma '|}\\
&& P_{\gamma'_1} (\partial^{(\gamma' _2)} (X)) (t,x,\omega)
\partial^{\gamma '_1} \varphi (X(t,x,\omega)) |^2 dx
\Eea where $P_{\gamma '_1} (x_1, \ldots x_d)$ is a polynomial as in
(3.3)and for a multi index $\alpha = (\alpha_1,\ldots \alpha_d) $,
we use the notation $\partial^{(\alpha)} (X) (t,x,\omega)
=(\partial^{\alpha _1} X_1 (t,x,\omega), \ldots
\partial^{\alpha_d}X_d (t,x,\omega))$. Using the change of variable
$ y = X(t,x,\omega)$, the integral in the last equality above is
\Bea &&\leq \sum\limits_{|\gamma |+|\gamma '|=|\beta |}
~~\sum\limits_{|\gamma '_1|, |\gamma '_2| \leq |\gamma '|}
\int\limits_{K'_t} |(X^{-1} (t,x,\omega))^\alpha
\partial^\gamma \phi' (X^{-1} (t,x,\omega)) \\
&&P_{\gamma '} (\partial^{(\gamma '_2)} (X) (t,X^{-1}
(t,x,\omega),\omega)) \partial^{\gamma '_1} \phi (x)|^2~ |det
(\partial X(t,\omega)^{-1})(X^{-1} (t,x,\omega))|~dx. \Eea

where $\partial X(t,\omega)$ is the Jacobian of $x \rightarrow
X(t,x,\omega)$ , $\partial X(t,\omega)^{-1}$ is the inverse of
$\partial X (t,\omega)$ and $K_t'(\omega) = X(t,K',\omega) $ is the
image of $K'$ under the map $X(t,.,\omega)$. Hence from the above,
we get \Bea \|x^\alpha \partial^\beta (\phi'
X_t(\omega)(\phi))\|^2_0 &\leq& C ~\alpha _1 (t) \alpha _2 (t) ~
\sum\limits_{|\gamma'| \leq|\beta |}\int
|\partial^{\gamma'}\phi(x)|^2~dx \Eea for some constant $C$
depending only on $p$ and $K'$ , where \Bea \alpha _1 (t,\omega) &=&
\max_{|\beta |\leq 2p} \max_{|\gamma '_1|,|\gamma '_2|\leq |\beta |}
\sup\limits_{x \in K'} |P_{\gamma '_1} (\partial^{(\gamma '_2)}
(X))(t,x,\omega)|^2 \\
\alpha _2 (t,\omega) &=& \sup\limits_{x \in K'} |(det (\partial
X)^{-1} (t,x,\omega))|. \Eea Summing over $\alpha$ and $\beta$ with
$|\alpha|+|\beta| \leq 2p$, we get \Bea \|\phi'
X_t(\omega)(\phi))\|^2_p~ \leq ~ C~\alpha _1 (t) \alpha _2 (t)
\|\phi\|^2_p \Eea Hence \Bea E\|\phi'X_t(\phi))\|_p^2 &\leq& C
\|\phi\|^2_p E(\alpha_1(t)
\alpha_2(t)) \\
&\leq& C \|\phi\|^2_p\sup\limits_{t \leq T} ~(E(\alpha _1
(t))^2)^{1/2}(E(\alpha _2
(t))^2)^{1/2} \\
&\leq& C\|\phi\|^2_p \Eea for some constant $C$ that changes from
line to line, but depends only on $p,T,$ and $K$. Note that we have
used Theorem 2.2 in the last inequality. It now follows that \Bea
|\langle S_t (\psi), \phi \rangle | \leq C\|\psi\|_{-p} \|\phi\|_p
\Eea and this completes the proof .

\hfill$\Box$

\end{document}